\documentclass[leqno]{amsart}
\usepackage{amsfonts,amssymb,amsmath,amsgen,amsthm}
\usepackage{hyperref}
\usepackage{color}

\theoremstyle{plain}
\newtheorem{theorem}{Theorem}[section]

\newtheorem{corollary}[theorem]{Corollary}
\newtheorem{proposition}[theorem]{Proposition}

\theoremstyle{remark}
\newtheorem{remark}[theorem]{Remark}

\newtheorem{example}[theorem]{Example}

\def\dis
{\displaystyle}

\def\R{{\mathbb R}}

\def\({\left(}
\def\){\right)}
\def\<{\left\langle}
\def\>{\right\rangle}
\def\le{\leqslant}
\def\ge{\geqslant}

\def\Tend#1#2{\mathop{\longrightarrow}\limits_{#1\rightarrow#2}}

\def\d{{\partial}}
\def\eps{\varepsilon}

\DeclareMathOperator{\IM}{Im}

\numberwithin{equation}{section}

\begin{document}

\title[Explicit solutions]{Explicit solutions for replicator-mutator
  equations: extinction vs. acceleration}

\author[M. Alfaro]{Matthieu Alfaro}
\address{CNRS \& Univ. Montpellier~2\\Math\'ematiques
\\CC~051\\34095 Montpellier\\ France}
\email{Matthieu.Alfaro@univ-montp2.fr}

\author[R. Carles]{R\'emi Carles}
\address{CNRS \& Univ. Montpellier~2\\Math\'ematiques
\\CC~051\\34095 Montpellier\\ France}
\email{Remi.Carles@math.cnrs.fr}

\begin{abstract} We consider a class of nonlocal reaction-diffusion
problems, referred to as {\it replicator-mutator} equations in
evolutionary genetics. By using explicit changes of unknown
function, we show that they are equivalent to the heat equation
and, therefore, compute their solution explicitly. Based on this,
we then prove that, in the case of beneficial mutations in asexual
populations, solutions dramatically depend on the tails of the
initial data: they can be global, become extinct in finite time
or, even, be defined for no positive time. In the former case, we
prove that solutions are accelerating, and in many cases converge for
large time to
some universal Gaussian profile. This sheds light on the
biological relevance of such models.
\end{abstract}



\maketitle


\section{Introduction}
\label{sec:introduction}

We consider {\it replicator-mutator} equations, that is nonlocal
reaction-diffusion problems of the form
\begin{equation}\label{eq-bio-poids}
\partial _t u=\partial _{xx}u+\left(f(x)-\int _\R f(x)u(t,x)\,dx\right)u, \quad t>0,\; x\in \R,
\end{equation}
where $f(x)$ is a given weight. In this context, $u(t,x)$ is the
density of a population (at time $t$ and per unit of fitness) on a
one-dimensional fitness space. We detail below the biological
background of such models.

\medskip

In this work, we mainly focus on \eqref{eq-bio-poids} for the
special case $f(x)=x$, namely
\begin{equation}\label{eq-bio}
\partial _t u=\partial _{xx}u+(x-\bar u(t))u, \quad t>0,\; x\in \R,
\end{equation}
where the nonlocal term is given by
\begin{equation}\label{premier-moment}
\bar u(t):=\int _\R xu(t,x)\,dx.
\end{equation}
We make a rigorous and detailed analysis of the Cauchy problem
associated with \eqref{eq-bio}. Precisely, we prove that it can be
reduced to the heat equation, and therefore compute its solution
explicitly. This enables us to describe a variety of contrasted
behaviors (extinction, acceleration\dots) depending on the initial
data.

\begin{remark}[Generalizations to quadratic weights]\label{rem:quadratic}  As a matter of fact, our analysis is also valid for
quadratic weights. Following the algebraic reductions of Section
\ref{sec:algebraic-reductions}, one may easily collect explicit
formulas for the solutions of \eqref{eq-bio-poids}  when $f(x)=\pm
x^2$, and, based on this, explore their behaviors. Nevertheless,
since the model \eqref{eq-bio} triggered a flow of studies in
evolutionary genetics, we just state our results
for this well-established case.
\end{remark}

In the context of evolutionary genetics, equation \eqref{eq-bio}
was introduced by Tsimring  et al. \cite{TLK96}, where they
propose a mean-field theory for the evolution of RNA virus
populations on a fitness space. Without mutations, and under the
constraint of constant mass $\int _\R u(t,x)\,dx=1$, the dynamics
is given by
$$
\partial _t u=(x-\bar u(t))u,
$$
where $\bar u(t)=\int _\R xu(t,x)\,dx$ is the average fitness of
the virus population. As a first step to take into account
evolutionary phenomena, one can then model mutations by the
Laplace diffusion operator so that the above integro-differential
equation is transferred into \eqref{eq-bio}. Notice also that
equation~\eqref{eq-bio} appears as a mean-field model for
diffusion-limited growth \cite{WS83}.

A central issue in evolutionary genetics is to predict whether a
population accumulates deleterious or advantageous mutations. The
former case is known as the {\it Muller's ratchet} \cite{Muller1,
Muller2}: an asexual population will accumulate deleterious
mutations and, therefore, its fitness will decay. On the other
hand, it recently turned out that beneficial mutations are more
abundant than previously suspected. Hence, after the seminal work
\cite{TLK96}, equation~\eqref{eq-bio} received a lot of attention
since it enables to capture the effect of such beneficial
mutations in asexual (clonal) populations. For more details and
comments on biological assumptions and such models, we refer to
\cite{RWC02,GCPS07,RBW08}, the review \cite{SG10}, \cite{WI13} and
the references therein.

However, for biological applications, the unlimited growth rate of
$u(t,x)$ at large $x$ in \eqref{eq-bio} is not admissible. To deal
with such a problem, the authors of the aforementioned works
consider a \lq\lq cut-off version'' of \eqref{eq-bio} at large
$u(t,x)$ \cite{TLK96,RWC02,SG10}, or provide a
proper stochastic treatment for large fitness region \cite{RBW08}.
In the former cut-off regime, the existence of solitary waves
(that is localized nonnegative profiles travelling at constant
speed and shape) and the way they attract solutions of the Cauchy
problem are investigated. In particular, the speed of the wave is
determined by a matching condition, and solutions of the Cauchy
problem travel at this constant speed in the large time regime.

\medskip

We now go back to the original deterministic equation
\eqref{eq-bio}. As far as we know, little was known concerning
existence and behaviors of solutions. Let us here mention the main
result of Biktashev \cite{B14}: for compactly supported initial
data, solutions converge, as $t\to \infty$, to a Gaussian profile,
where the convergence is understood in terms of the moments of
$u(t,x)$. One may then conjecture that this property remains valid
for \lq\lq arbitrary'' initial data. In this work, we show in
particular that this is completely false: tails of the initial
data have a strong influence on solutions.

The situation for equation \eqref{eq-bio} is also in sharp
contrast with the cut-off and stochastic approximations as studied
in \cite{TLK96,RWC02,SG10,RBW08}. First, using the Fourier
transform, one can explicitly compute all solitary waves and
observe that not only all positive speeds are admissible but also
that all profiles are changing sign (see
Appendix~\ref{sec:solit-waves-eqref} for details). Next, solutions
of the Cauchy problem can become extinct in finite time and, if
global, are accelerating as time passes. This is the main goal of
this work to rigorously prove these features for \eqref{eq-bio}.

\medskip

Throughout this work, we assume that the initial data is nonnegative,
$u_0(x)\ge 0$, and satisfies
\begin{equation}
\int _\R u_0(x)\,dx=1,
\end{equation}
so that, \emph{formally}, $\int _\R u(t,x)\,dx=1$ for later times. Indeed, if we
formally integrate \eqref{eq-bio} over $x\in \R$, we see
that the total mass $m(t):=\int _\R u(t,x)\,dt$ solves the Cauchy
problem
\begin{equation}\label{ode-masse}
\frac d{dt}m(t)=(1-m(t))\bar u (t), \quad m(0)=1,
\end{equation}
so that the Gronwall lemma yields $m(t)=1$ as long as $\bar u(t)$ is
meaningful. A striking result of this paper is that the above formal
argument may turn out to be completely wrong, in the sense that the
solution may become extinct in finite time, $u(t,x)=0$ for all $x\in
\R$ and $t\ge T$.

\medskip

The organization of the paper is as follows. In Section
\ref{sec:results}, we state our main results for \eqref{eq-bio}.
The keystone result is Theorem~\ref{th:explicit-sol} and contains
explicit formulas for solutions. Its proof (and that of some
generalizations as explained in Remark~\ref{rem:quadratic})
involves algebraic reductions that are given in
Section~\ref{sec:algebraic-reductions}. The different scenarii for
solutions (extinction in finite time, global existence,
acceleration\dots) are then proved in
Section~\ref{sec:cauchy-problem}. We give a short summary of our
work in Section \ref{sec:conclusion}. Last, the solitary waves are
computed in Appendix~\ref{sec:solit-waves-eqref}, and the
propagation of Gaussian initial data in the case of a quadratic
weight $f$ in \eqref{eq-bio-poids} is presented in
Appendix~\ref{sec:prop-gauss-init-quad}.

\section{Main results}
\label{sec:results}

By using tricky algebraic manipulations, we can actually reduce
the nonlocal equation \eqref{eq-bio} to the heat equation, and
therefore compute the solution explicitly. 
This is our first main
result and it reads as follows.

\begin{theorem}[The solution explicitly]
\label{th:explicit-sol} Let $u_0\ge 0$, with $\int u_0=1$. As long
as $\bar u(t)$ is finite, the solution of \eqref{eq-bio} with
initial data $u_0$ is given by
\begin{equation}\label{formula-avec-w}
u(t,x)= \frac{e^{tx+t^3/3}w(t,x+t^2)}{1+\displaystyle\int
_0^t \int _\R xe^{sx+s^3/3}w(s,x+s^2) \,dx\,ds},
\end{equation}
where $w(t,x)=e^{t\d_{xx}}u_0(x)$ is the solution of the heat equation with initial data
$u_0$. As a consequence, we also have
\begin{equation}\label{formula-grosse}
u(t,x)= \frac{e^{tx+t^3/3} \displaystyle\int _\R \frac 1
{\sqrt{4\pi t}}e^{-(x+t^2-y)^2/4t}u_0(y)\,dy}{1+\displaystyle\int
_0^t \int _\R xe^{sx+ s^3/3}\int _\R \frac 1 {\sqrt{4\pi
s}}e^{-(x+s^2-y)^2/4s}u_0(y)\,dy \,dx\,ds},
\end{equation}
and
\begin{equation}\label{formula-belle}
 u(t,x)=\frac{e^{tx}\displaystyle \int _\R
\frac 1 {\sqrt{4\pi
t}}e^{-(x+t^2-y)^2/(4t)}u_0(y)\,dy}{\displaystyle \int _\R
e^{ty}u_0(y)\,dy}.
\end{equation}
\end{theorem}

\begin{corollary}[The nonlocal term explicitly]\label{cor:u-bar} As long as it exists, $\bar
u(t)$ is given by
\begin{equation}\label{u-bar}
\bar u(t)=t^2+\frac{\displaystyle \int _\R
e^{ty}y\, u_0(y)\,dy}{\displaystyle \int _\R e^{ty}u_0(y)\,dy}.
\end{equation}
\end{corollary}
It seems however that these explicit formulas rely on the fact that
the equation has \emph{exactly} the form \eqref{eq-bio}: if saturation
(as in \cite{TLK96,RWC02,SG10}) or stochasticity (as in \cite{RBW08})
is introduced, then we can no longer take advantage of this
``algebraic miracle''.  
\smallbreak

Equipped with the above formulas, we can prove rather different
scenarii for the Cauchy problem associated with \eqref{eq-bio}.
Let us notice that, without our exact formulas, proving such
behaviors seems to be far from obvious.
\begin{theorem}[Global existence vs. extinction in finite time]\label{th:global}
Let $u_0\ge 0$, with $\int u_0=1$. Consider
\begin{equation*}
  T = \sup\left\{ t\ge 0,\quad \int_0^\infty e^{ty}u_0(y)dy<\infty\right\}\in [0,\infty].
\end{equation*}
\begin{itemize}
\item [$(i)$]  If $T=\infty$, then in \eqref{eq-bio}, both
$u(t,x)$ and $\bar u(t)$ are global in time. Typically, $u \in
L^\infty_{\rm loc}((0,\infty)\times \R)$, $\bar u\in L^\infty_{\rm
loc}(0,\infty)$, and $\int_\R u(t,x)dx=1$ for all $t\ge 0$. \item
[$(ii)$] If $0<T<\infty$ ,  then extinction in finite time occurs,
that is
\begin{equation*}
  u(t,x)=0,\quad \forall t> T, \ \forall x\in \R.
\end{equation*}
\item [$(iii)$] If $T=0$, then $u(t,x)$ is defined for no $t>0$.
\end{itemize}
\end{theorem}

The first case holds, for instance, for Gaussian initial data
whose propagation is investigated in Proposition
\ref{prop:prop-gauss}. The proof of $(i)$ is obvious since the
assumption $\int _0 ^\infty e^{ty} u_0(y)\,dy <\infty$ for all
$t>0$ (i.e. $T=\infty$) implies $\int _0 ^\infty e^{ty}y
u_0(y)\,dy <\infty$ for all $t>0$, and therefore (notice that the
integration on $(-\infty,0)$ is harmless since $ye^{ty}$ is
bounded on this interval) both \eqref{u-bar} and
\eqref{formula-belle} are meaningful for all $t>0$.

On the other hand, initial data not having {\it very light} tails
at $+\infty$ make the equation completely meaningless in positive
and finite time (second point). This in particular happens for
initial data having {\it light} exponential tails. The proof of
$(ii)$ is straightforward in view of \eqref{formula-belle}:
\begin{equation*}
  0\le u(t,x) \le \frac{e^{tx}}{\sqrt{4\pi t}}\frac{\dis\int_\R u_0(y)dy}{\displaystyle \int _\R
e^{ty}u_0(y)\,dy}= \frac{e^{tx}}{\sqrt{4\pi t}\displaystyle \int _\R
e^{ty}u_0(y)\,dy}.
\end{equation*}
The numerator
remains bounded for each $(t,x)$ fixed, while the denominator tends to
$+\infty$ near time $T$ (possibly just after $T$).
\begin{example}[Light exponential tail, extinction for $t\ge \alpha$]\label{ex:exp}
  If
\begin{equation}\label{exponential-data}
u_0(y)=\alpha e^{-\alpha y} \mathbf 1 _{(0,\infty)}(y),\quad
\alpha >0,
\end{equation}
then in \eqref{eq-bio},  both $u(t,x)$ and $\bar u(t)$ are defined
on $(0,\alpha)$. They are given (see subsection
\ref{subsec:extinction} for details) by
\begin{equation}\label{u-bar-extinction}
\bar u(t)=t^2+\frac 1{\alpha -t} \Tend t \alpha  \infty,
\end{equation}
and
\begin{equation}\label{u-extinction}
u(t,x)=\frac 1{\sqrt{2\pi}}(\alpha -t)e^{-(\alpha -t)x}e^{-\alpha
t^2+\alpha ^2 t} {\rm Erf}\left(\frac{-(x+t^2-2\alpha
t)}{\sqrt{2t}}\right)  \Tend t \alpha  0,
\end{equation}
uniformly in $x\in \R$, where
\begin{equation*}
{\rm Erf} (\theta):=\int _\theta ^\infty
e^{-z^2/2}\,dz.
\end{equation*}
 In view of this, it seems reasonable to extend the solution by
$u(t,x)\equiv 0$ for $t\ge \alpha$, which shows an extinction
phenomena.
\end{example}
\begin{example}[Light tail, extinction for $t>\alpha$]
  Consider a slight modification of the above example:
\begin{equation*}
u_0(y)=\tilde\alpha \frac{1}{y^2} e^{-\alpha y} \mathbf 1 _{(1,\infty)}(y),\quad
\alpha >0,
\end{equation*}
where $\tilde \alpha$ is chosen so that $\int_{\R} u_0=1$.
Invoking \eqref{formula-belle}-\eqref{u-bar}, the formula that we
obtain is not as explicit as  \eqref{u-extinction}. However, it is
clear that $\bar u(t)$ is finite for $t\le \alpha$, while
$u(t,x)\equiv 0$ for $t>\alpha$.
\end{example}

Last, as suggested by the denominator of
formula \eqref{formula-belle}, initial data having {\it heavy}
tails prevent the definition of the solution for any positive
time, that is $(iii)$. See Remarks~\ref{rem:T=0} and \ref{rem:T=0bis}
for a precise explanation.

\begin{example}[Heavy tails]\label{extinction-immediate}
If
\begin{equation}\label{hyp-immediate}
y\mapsto e^{ty}  u_0(y) \notin L^1(0,\infty), \quad \forall t>0,
\end{equation}
then the solution $u(t,x)$ of \eqref{eq-bio} is defined for no
$t>0$. This is typically the case if $u_0$ decays only
algebraically.
\end{example}

\begin{remark}
  The fact that not enough decay of the initial data on one side leads
  to pathological phenomena can be compared to a situation recently
  studied in the framework of dispersive equations. For the
  $L^2$-critical generalized Korteweg-de Vries equation
  \begin{equation*}
    \d_t u + \d_x \(\d_{xx}u +u^5\)=0;\quad ;\quad u_{\mid t=0}=Q+\eps_0,
  \end{equation*}
where $Q$ is the unique even positive solution to $Q''+Q^5=Q$, given by
\begin{equation*}
  Q(x)=\( \frac{3}{\cosh^2(2x)}\)^{1/4},
\end{equation*}
Martel, Merle and Rapha\"el \cite{MaMeRa-p} have proved that if the
initial perturbation $\eps_0$ does not decay sufficiently fast on the right, then various
 regimes are possible, including a continuum of blow-up rates, a
 continuum of growth rate at infinity, while if $\int_0^\infty
x^{10}\eps_0(x)dx<\infty$, then only three scenarii are possible. In
the case of the parabolic energy critical harmonic heat flow, similar
phenomena had been observed by Gustafson, Nakanishi and Tsai \cite{GuNaTs10}.
\end{remark}
Let us now turn to the speed of propagation of solutions. Plugging
$u_0(y)=\delta _0(y)$ the Dirac mass at $0$ in
\eqref{formula-belle} and \eqref{u-bar}, one gets (see subsection
\ref{subsec:long-time} for details)
\begin{equation}\label{sol-elementaire}
u(t,x)=\frac 1{\sqrt{4\pi t}}e^{-(x-t^2)^2/4t}, \quad \bar
u(t)=t^2.
\end{equation}
This suggests that the solution of the Cauchy problem are
accelerating. To maintain this affirmation, we investigate the
propagation of a Gaussian initial data, which is relevant for
biological lectures.
\begin{proposition}[Accelerating propagation of Gaussian initial
data]\label{prop:prop-gauss} If
\begin{equation}\label{initial-gauss}
u_0(x)=\sqrt{\frac a{2\pi}}e^{-a(x-m)^2/2},\quad a>0, \quad m\in
\R,
\end{equation}
then the solution of \eqref{eq-bio} is
\begin{equation}\label{sol-gauss-cond-ini}
u(t,x)=\sqrt{\frac {a(t)}{2\pi}}e^{-a(t)(x-m(t))^2/2}, \quad
a(t):=\frac a{1+2at}, \quad m(t):=m+t^2+\frac t a.
\end{equation}
\end{proposition}
This shows that, starting from a Gaussian profile, the solution
remains a Gaussian function, is accelerating and flattening since
$m(t)\sim t^2$, $a(t)\sim \frac 1{2t}$, as $t\to \infty$.
Starting from our explicit formula, the computations that prove
the above proposition are presented in subsection
\ref{subsec:calcul-proposition}. Notice that this family of
Gaussian self-similar solutions already appears in \cite{B14},
where the long time convergence of the solution of \eqref{eq-bio}
(with a compactly supported initial data) to a Gaussian profile is
also investigated. As far as this result is concerned, we can
provide a sharp improvement of the convergence procedure.
Precisely, the long time convergence in \cite[Theorem~1]{B14} is
understood in term of the moments of $u(t,x)$, whereas we can
prove strong uniform convergence.
Precisely the following holds.

\begin{theorem}[Long time behavior for compactly supported initial data]\label{th:compactly-supp}
Let $u_0\ge 0$ be compactly supported, with $\int u_0=1$. Let
$u(t,x)$ be the global solution of \eqref{eq-bio} with initial
data $u_0$. Then there is $C>0$ such that
$$
\sup _{x\in \R}\; \left |u(t,x)-\frac 1{\sqrt{4\pi
t}}e^{-(x-t^2)^2/4t} \right|\le \frac C t,\quad \forall t\ge 1.
$$
\end{theorem}

The above result actually measures, uniformly with respect to $x\in \R$,
the deviation from the elementary solution
\eqref{sol-elementaire}. The proof is based on a combination of
our explicit formulas with an elementary estimate on the long time
behavior of the heat equation. It will appear in subsection
\ref{subsec:long-time}.

\section{Algebraic reductions}
\label{sec:algebraic-reductions}

In this section, we show how to relate the solution of various
modulations of \eqref{eq-bio} with the solution of the standard
heat equation
\begin{equation}
  \label{eq:heat}
  \d_t w = \d_{xx}w, \quad t>0,\; x\in \R;\quad w_{\mid t=0}=u_0,
\end{equation}
or a perturbation of the heat equation. In particular, the proof
of the main result Theorem \ref{th:explicit-sol} will appear in
subsection \ref{subsec:spat-line-fact}.

\subsection{External time-dependent factor}
\label{subsec:extern-time-depend}

Consider the equation
\begin{equation}
  \label{eq:ext}
  \d_t u = \d_{xx}u +a(t)u+g(t,x)u, \quad t>0,\; x\in \R;\quad u_{\mid t=0}=u_0,
\end{equation}
where $a$ is a given function of time only (independent of $x$ and
$v$), and $g$ is independent of $u$. Consider $v$ the solution to the
Cauchy problem
\begin{equation}
  \label{eq:ext2}
  \d_t v = \d_{xx}v +g(t,x)v, \quad t>0,\; x\in \R;\quad v_{\mid t=0}=u_0.
\end{equation}
Then $u$ and $v$ are explicitly related through the formula
\begin{equation*}
  u(t,x) = v(t,x)e^{\int_0^t a(s)ds}.
\end{equation*}

\subsection{Generalized momentum factor}
\label{subsec:gener-moment-fact}

Suppose now that in \eqref{eq:ext}, the time dependent function is
related to $u$ in the same fashion as in \eqref{premier-moment},
\begin{equation}\label{eq:moment-gen}
  \d_tu = \d_{xx}u+g(t,x)u-\overline u(t)u, \quad t>0,\; x\in \R;\quad u_{\mid t=0}=u_0,
\end{equation}
where
\begin{equation*}
  \overline u(t)=\int_\R f(x)u(t,x)dx,
\end{equation*}
for some weight function $f(x)$. Introduce $v$ the solution to the
Cauchy problem
\begin{equation}\label{eq:moment-gen2}
  \d_tv = \d_{xx}v+g(t,x)v, \quad t>0,\; x\in \R;\quad v_{\mid t=0}=u_0.
\end{equation}
 Then formally,
\begin{equation*}
  v(t,x) = u(t,x)e^{\int_0^t \overline u(s)ds}.
\end{equation*}
We remark that this change of unknown function can be inverted:
multiplying the above expression by $f(x)$ and integrating over $x\in
\R$, we get
\begin{equation*}
  \overline v(t) = \overline u(t) e^{\int_0^t \overline u(s)ds} =
  \frac{d}{dt}\( e^{\int_0^t \overline u(s)ds}\) .
\end{equation*}
By integrating in time, we infer
\begin{equation*}
  \int_0^t \overline v(s)ds =   e^{\int_0^t \overline u(s)ds} -1,
\end{equation*}
and, so long as $\int_0^t \overline v(s)ds>-1$,
\begin{equation}\label{eq:chgt1}
  u(t,x) =\frac{v(t,x)}{1+\dis \int_0^t \overline v(s)ds}.
\end{equation}
In the case considered throughout this paper, $u_0\ge 0$, which
implies, as we will see below,
$v(t,x)>0$ for all $t>0$ and all $x\in \R$ in the case
$g(t,x)=x$. Therefore, we always have
$\int_0^t \overline v(s)ds\ge 0$, and the above computations are licit
provided that $\overline u$ (and therefore $\overline v$) is finite.
\begin{example}\label{ex1}
  Consider \eqref{eq-bio} without the drift factor $xu$, that is
  \begin{equation*}
     \d_t u = \d_{xx}u -\overline u(t)u, \quad t>0,\; x\in \R;\quad v_{\mid t=0}=u_0,
  \end{equation*}
with
\begin{equation*}
  \overline u(t)=\int_\R x u(t,x)dx.
\end{equation*}
In that case, $v=w$, solution to the heat equation \eqref{eq:heat}.
In view of the expression of the heat kernel, we have:
\begin{equation*}
  \int_0^t \overline w(s)ds = \int_0^t \int _\R \int _\R x\frac{1}{\sqrt{4\pi s}}
  e^{-(x-y)^2/4s}u_0(y)\,dy\,dx\,ds.
\end{equation*}
We compute
\begin{align*}
  \int_{\R}x  e^{-(x-y)^2/4s} dx  = \int_{\R}(x+y)  e^{-x^2/4s} dx =
  y\sqrt{4\pi s},
\end{align*}
and thus
\begin{equation}\label{bidule}
  u(t,x) = \frac{w(t,x)}{1 +t \dis \int_{\R} y u_0(y)dy}.
\end{equation}
Therefore, if $u_0$ is even, or if the main part of its mass lies on
the right, $\int_\R y u_0(y)dy\ge 0$, then the solution $u$ is
well-defined for all times $t\ge 0$. On the other hand, if the mass of
$u_0$ is more important on the left, $\int_\R y u_0(y)dy< 0$, then
finite time blow-up occurs:
\begin{equation*}
  \exists T_*>0,\quad u(t,x)\Tend t {T_*} +\infty,\quad \forall x\in\R.
\end{equation*}
\end{example}

\begin{remark}\label{rem:T=0}
  The above reduction requires to be able to consider an open time
  interval, in order for the integration procedure to make sense. This
  approach becomes meaningless if we have $\overline u(t)=\infty$
  (hence $\overline v(t)=\infty$) for all $t>0$, which is exactly the
  case of Theorem~\ref{th:global}, {\it (iii)}.
\end{remark}

\subsection{Spatially linear factor}
\label{subsec:spat-line-fact}

Consider now a heat equation supplemented with an extra term involving
a factor which is linear in $x$,
\begin{equation}
  \label{eq:ext-lin}
  \d_t v = \d_{xx}v +a(t)x v, \quad t>0,\; x\in \R;\quad v_{\mid t=0}=u_0,
\end{equation}
where $a$ is a given function of time only (independent of $x$ and
$v$). In quantum mechanics, the left hand side of the equation is replaced
by $i\d_t v$, where $i=\sqrt{-1}$,  and the corresponding Schr\"odinger
equation models the evolution of particles under the effect of an electric field
$a(t)x$. When the function $a$ is constant, it
is possible to relate the solution of the free Schr\"odinger equation
to the solution of the equation with this electric field through the
\emph{Avron--Herbst formula}, see e.g. \cite{Thirring}. This formula
can be generalized to the case where $a$ does depend on $t$, see
\cite{Ca11}. Replacing $t$ with $-it$ in the formula given in
\cite{Ca11}, we see that the solutions to \eqref{eq:heat} and
\eqref{eq:ext-lin} are related through
\begin{equation*}
   v(t,x) = w\(t,x + 2 \int_0^t\int_0^s a(\tau)d\tau ds\)\exp\(
  x\int_0^t a(s)ds +\int_0^t \(\int_0^s a(\tau)d\tau\)^2  ds\).
\end{equation*}
In the case $a(t)\equiv 1$, this formula is simply
\begin{equation}
  \label{eq:AH-gen}
   v(t,x) = w\(t,x + t^2\)\exp\(
  t x+\frac{t^3}{3}\).
\end{equation}

\begin{proof}[Proofs of Theorem~\ref{th:explicit-sol} and Corollary \ref{cor:u-bar}]
Combining \eqref{eq:chgt1} and \eqref{eq:AH-gen}, we infer
\eqref{formula-avec-w}. The expression \eqref{formula-grosse} then
stems from the explicit formula of the heat kernel on $\R$.
Finally, to deduce \eqref{formula-belle}, we denote by $I(t)$ the
triple integral appearing in the denominator of
\eqref{formula-grosse}. Using Fubini's Theorem, we first compute
the integral with respect to $x$. Using elementary algebra
(canonical form) we find
\begin{align}
\int _\R xe^{sx}e^{-(x+s^2-y)^2/4s}\,dx&= e^{sy}\int _\R
xe^{-[x-(s^2+y)]^2/(4s)}\,dx \nonumber\\
&= e^{sy}\int _\R (z+(s^2+y))e^{-z^2/(4s)}\,dz\nonumber \\
&= e^{sy}(s^2+y)\sqrt{4\pi s}.\label{petit-calcul}
\end{align}
As a result, we have
$$
I(t)=\int _0 ^t \int _\R (s^2+y)e^{
s^3/3+sy}u_0(y)\,dy\,ds=\int _\R (e^{t^3/3+ty}-1)u_0(y)\,dy.
$$
Plugging this into \eqref{formula-grosse} and using the normalization $\int_\R
u_0=1$, we then obtain \eqref{formula-belle}. Using
\eqref{formula-belle} and equality \eqref{petit-calcul}
again, we see that \eqref{u-bar} holds true.
\end{proof}

\begin{remark}
  The denominator of \eqref{bidule} in Example~\ref{ex1} corresponds
  to the expression
  obtained by considering the first two terms of the Taylor expansion
  of the exponential in the denominator in
  \eqref{formula-belle}. Example~\ref{ex1} illustrates the fact
  that introducing the term $xu$ in \eqref{eq-bio} prevents
  blow-up, as shown by the formula \eqref{formula-belle}  and
  Theorem~\ref{th:global}.
\end{remark}

\begin{remark}\label{rem:T=0bis}
  Back to Theorem~\ref{th:global}, \emph{(iii)}, we see that if there
  was a $\tau>0$ such that $\overline u$ is finite on $[0,\tau]$, then
  \eqref{formula-belle} would hold true. On the other hand, the
  assumption $T=0$, along with \eqref{formula-belle}, would imply
  $u(t,x)=0$ for all $t\in (0,\tau]$ and all $x\in \R$, while we have
  seen in \eqref{ode-masse} that so long as $\overline u$ is finite, we have
  $\int_\R u(t,x)dx=1$, hence a contradiction.
\end{remark}
\subsection{Spatially quadratic factor}
\label{subsec:spat-quadr-fact} Consider
\begin{equation}
  \label{eq:quad}
  \d_t v= \d_x^2 v -a(t)x^2v\quad;\quad v_{\mid t=0}=u_0,
\end{equation}
where $a$ is a given function of time only (independent of $x$ and
$v$). In the case where $a$ is constant (say $a=1$), the solution to
\eqref{eq:quad} is given by the \emph{Mehler's formula},
\begin{equation}
  \label{eq:mehler}
  v(t,x) = \frac{1}{\sqrt{2\pi \sinh(2t)}}\int_\R e^{-\coth (2t)
    \frac{x^2+y^2}{2}-\mathrm{cosech} (2t)xy} u_0(y)dy.
\end{equation}
The formula is known in the context of the heat equation (\cite{Gr09})
as well as in the context of the Schr\"odinger equation
(\cite{Feyn}). For a general time-dependent function $a$, introduce
the fundamental solution associated to the corresponding oscillator,
\begin{equation*}
  \label{eq:solfond}
\left\{
  \begin{aligned}
&\ddot \mu -a(t)\mu =0 \quad;\quad \mu(0)=0,\quad
  \dot\mu(0)=1, \\
 &   \ddot \nu - a(t)\nu =0 \quad;\quad \nu(0)=1,\quad
  \dot\nu(0)=0.
  \end{aligned}
  \right.
\end{equation*}
For $a(t)\ge 0$, we check that $\nu(t)\ge 1$ for all $t\ge 0$, and
$\mu(t)>0$ for all $t>0$. Adapting the generalized lens transform
presented in \cite{Ca11}, we
see that the solutions to \eqref{eq:quad} and \eqref{eq:heat} are
related through the formula
\begin{equation}\label{eq:lens}
  v(t,x) = \frac{1}{\sqrt{\nu(2t)}} e^{-\frac{x^2}{2}\frac{\dot \nu(2t)}{\nu(2t)}}
  w\(\frac{\mu(2t)}{2\nu(2t)},\frac{x}{\nu(2t)}\).
\end{equation}
Of course, this formula makes sense so long as $\nu$ is nonzero, and
so long as the map $t\mapsto \mu(2t)/\nu(2t)$ is invertible. Note that
this is the case for all positive times when $a\ge 0$, from the above
remark.
\begin{remark}
  In the case $a=1$, we compute explicitly $\mu(t) = \sinh(t)$
and $\nu(t)=\cosh(t)$. Mehler's formula \eqref{eq:mehler} can be
viewed as the composition of the lens transform \eqref{eq:lens} and
the explicit formula for the heat kernel.
\end{remark}

\begin{remark}[Multidimensional case]
  All the formulas presented in this section can be generalized to a
  multidimensional framework, $x\in \R^d$, $d\ge 1$. In the case
  considered in subsection \ref{subsec:spat-line-fact}, replace $a(t)x$ with
  $a(t)\cdot x$ where $a(t)\in \R^d$ is a vector-valued
  time-dependent function.  In the quadratic case of
  subsection \ref{subsec:spat-quadr-fact}, it seems necessary to restrict to
  the isotropic case where $a(t)x^2u $ is replaced by
  \begin{equation*}
    a(t)|x|^2 u=a(t)\(\sum_{j=1}^d x_j^2\)u,
  \end{equation*}
that is, the coefficient in factor on $x_j^2$ is independent of $j$
(see \cite{Ca11}).
\end{remark}

\section{Proofs of various features of the Cauchy problem}
\label{sec:cauchy-problem}

In this section, based on our explicit formulas, we prove the
different behaviors as stated in Section \ref{sec:results}.

\subsection{Extinction in finite time}
\label{subsec:extinction}

We present here the computations associated to Example~\ref{ex:exp}. For the
initial data \eqref{exponential-data}, we compute
$$
\int _\R e^{ty}u_0(y)\,dy=\frac \alpha{\alpha -t}, \quad \int _\R
e^{ty}yu_0(y)\,dy=\frac \alpha{(\alpha -t)^2},
$$
which we plug into \eqref{u-bar} to get \eqref{u-bar-extinction}.
Next, \eqref{formula-belle} and elementary algebra (canonical
form) yields
\begin{eqnarray*}
 u(t,x)&=&
 \frac{e^{tx}(\alpha -t)}{\sqrt{4\pi t}}\int _0^\infty
 e^{-(x+t^2-y)^2/(4t)}e^{-\alpha y}\,dy\\
 &=& \frac{e^{tx}(\alpha -t)}{\sqrt{4\pi t}}
 e^{-\alpha(x+t^2-\alpha t)}\int _0^\infty e^{-[y-(x+t^2-2\alpha
 t)]^2/(4t)}\,dy\\
 &=& \frac{e^{tx}(\alpha -t)}{\sqrt{2\pi }}e^{-\alpha(x+t^2-\alpha t)} \int
 _{-\frac{x+t^2-2\alpha t}{\sqrt{2t}}} ^\infty e^{-z^2/2}\,dz,
\end{eqnarray*}
that is formula \eqref{u-extinction}. The fact that $u(t,x)\to 0$,
as $t\to \alpha$, uniformly in $x\in \R$ follows from the
following two facts: first, if $x\ge -\frac 1 {\alpha -t}$ then
\eqref{u-extinction} implies $|u(t,x)|\le C (\alpha -t)$; next,
for $t$ sufficiently close to $\alpha$, if $x\le -\frac 1 {\alpha
-t}$ then, using ${\rm Erf}(\theta) \sim \frac 1 \theta e^{-\theta
^2/2}$ as $\theta \to \infty$, \eqref{u-extinction} implies that
$$
|u(t,x)|\le Ce^{-(\alpha -t)x}{\rm Erf}\left(- \frac x
{\sqrt{2\alpha}}\right)\le Ce^{-(\alpha -t)x} 2
\frac{\sqrt{2\alpha}}{-x}e^{-x^2/(4\alpha)},
$$
so that $|u(t,x)|\le C' e^{-\frac {x^2}{4\alpha}\left(1+4\alpha
\frac{\alpha -t}x\right)}\le C'e^{-\frac {x^2}{8\alpha}}\le
C'e^{-\frac 1{8\alpha(\alpha -t)^2}}$.

\qed

\subsection{Propagation of Gaussian initial data.}
\label{subsec:calcul-proposition}

We now present the straightforward computations that prove
Proposition \ref{prop:prop-gauss}. We plug the initial data
\eqref{initial-gauss} into the formula \eqref{formula-belle} for
$u(t,x)$ and denote by $N(t,x)$, $D(t)$ the numerator, denominator
respectively. Using elementary algebra (canonical form), we get
\begin{align*}
D(t)&=\sqrt{\frac a{2\pi}}\int _\R
e^{ty}e^{-a(y-m)^2/2}\,dy\\
&=\sqrt{\frac a{2\pi}}\int _\R e^{- a
[y-(m+\frac t a)]^2/2}e^{mt+t^2/(2a)}\,dy\\
&=e^{mt+t^2/(2a)},
\end{align*}
and
\begin{align*}
N(t,x)&=\frac 1 {\sqrt{4\pi t}} e^{tx}\sqrt{\frac a{2\pi}}
 \int _\R e^{-(x+t^2-y)^2/(4t)}e^{-a(y-m)^2/2}\,dy\\
 &=\frac 1 {\sqrt{4\pi t}} e^{tx}\sqrt{\frac a{2\pi}}
 \int _\R
 e^{-\frac{1+2at}{(4t)}(y-\frac{x+t^2+2amt}{1+2at})^2}e^{-\frac a2
 \frac{(x+t^2-m)^2}{1+2at}}\,dy\\
 &= \frac 1 {\sqrt{4\pi t}} e^{tx}\sqrt{\frac a{2\pi}}e^{-\frac a2
 \frac{(x+t^2-m)^2}{1+2at}}\sqrt{2\pi\frac{2t}{1+2at}}\\
 &=\sqrt{\frac a{2\pi}} \frac 1{\sqrt{1+2at}}e^{-\frac
 a{2(1+2at)}\left(x-\frac{at^2+am+t}a\right)^2}e^{\frac t 2
 \frac{2am+t}a}.
\end{align*}
Finally, $u(t,x)=N(t,x)/D(t)$ easily yields
\eqref{sol-gauss-cond-ini}. \qed

\subsection{Long time behavior for compactly supported initial data}\label{subsec:long-time}
We  now prove Theorem \ref{th:compactly-supp}. Using elementary
algebra, \eqref{formula-belle} is recast as
\begin{equation}\label{formula-belle-bis}
 u(t,x)=\frac1{\sqrt{4\pi
t}} \frac{\displaystyle \int _\R
e^{-(x-t^2-y)^2/(4t)}e^{ty}u_0(y)\,dy}{\displaystyle \int _\R
e^{ty}u_0(y)\,dy},
\end{equation}
which, in particular, is the appropriate form to compute the
elementary solution \eqref{sol-elementaire} for $u_0(y)=\delta _0
(y)$. It follows from \eqref{formula-belle-bis} that the deviation
$$\psi(t,x):=u(t,x)-\frac 1{\sqrt{4\pi t}}e^{-(x-t^2)^2/4t}$$
 is
given by
\begin{equation*}
\psi(t,x)= \frac{\displaystyle \int _\R
\(e^{-(x-t^2-y)^2/(4t)}-e^{-(x-t^2)^2/(4t)}\)e^{ty}u_0(y)\,dy}{\sqrt{4\pi
t}\displaystyle \int _\R e^{ty}u_0(y)\,dy}.
\end{equation*}
Using Taylor formula we write
\begin{align*}
\left|e^{-(x-t^2-y)^2/(4t)}-e^{-(x-t^2)^2/(4t)}\right|&=\left|\int _0 ^1
\frac {y}{\sqrt t} \frac{x-t^2-\theta y}{2\sqrt
t}e^{-(\frac{x-t^2-\theta y}{2\sqrt t})^2}\,d\theta\right|\\
&\le  \frac {|y|}{\sqrt t} \sup _{z\in \R} |ze^{-z^2}|,
\end{align*}
and get
$$
|\psi(t,x)|\le \frac C t \frac {\displaystyle \int _\R e^{ty}|y|
u_0(y)\,dy}{\displaystyle \int _\R e^{ty}u_0(y)\,dy}\le \frac C t
M,
$$
where $\text{supp}\, u_0 \subset [-M,M]$. Theorem
\ref{th:compactly-supp} is proved. \qed

\section{Brief summary}\label{sec:conclusion}

We are concerned with evolutionary genetics models for asexual
populations (viruses, microbes). In contrast with Muller's ratchet
we aim at understanding the dynamics when accumulation of
deleterious mutations is neglected. In order to incorporate the
effects of mutations, we use the nonlocal reaction-diffusion
deterministic model proposed in \cite{TLK96}, and referred to as
the replicator-mutator equation.

Our mathematical analysis shows that one can reduce the
replicator-mutator equation to the heat equation. As a result,
solutions are completely explicit which enables to prove various
nontrivial behaviors. First, for initial data with heavy tails,
the equation is immediately meaningless. Next, for light initial
tails, the solution becomes extinct in finite time, which violates
the mass constraint formally observed. Last, for very light
initial tails, we prove that solutions are global and are
accelerating as time passes. This prevents the convergence to a
solitary wave, as observed for some perturbations (cut-off
approximation or stochastic treatment) of the original equation.

\appendix

\section{Solitary waves for \eqref{eq-bio}}
\label{sec:solit-waves-eqref}

In this Appendix, we compute explicitly the solitary waves for
\eqref{eq-bio}. In particular, all positive speeds are admissible
and, the Airy function being involved, all solitary waves are
changing sign, which enforces some cut-off arguments for
applications to biology.

We plug the ansatz $u(t,x)=\phi(x-ct)$ into equation
\eqref{eq-bio}. We are therefore looking for a speed $c$ and a
profile $\phi$ such that
\begin{equation}\label{pb-soliton}
 \left\{\begin{array}{ll}
\phi ''+c\phi '+(x-\bar \phi )\phi=0\quad \text{ on } \R, \vspace{4pt}\\
\phi (\pm \infty)=0,\quad \displaystyle \int_ \R \phi=1,
         \end{array}
\right.
\end{equation}
where $\bar \phi:=\displaystyle \int _\R x\phi (x)\,dx$.

If $\phi$ solves \eqref{pb-soliton} then
\begin{equation}\label{def-psi}
\psi(x):=\phi(x+\bar \phi)
\end{equation}
solves
\begin{equation*}\label{pb-soliton-moment-nul}
 \left\{\begin{array}{ll}
\psi ''+c\psi '+x\psi=0\quad \text{ on } \R \vspace{4pt}\\
\psi (\pm \infty)=0,\quad \displaystyle \int_ \R \psi=1,\quad \bar
\psi=0.
         \end{array}
\right.
\end{equation*}
Applying Fourier transform to this linear problem yields
$$
-\xi^2 \hat \psi+ci\xi \hat \psi+i\frac{d\hat \psi}{d\xi}=0,\quad
\hat \psi (0)=1,
$$
which is solved as
$$
\hat \psi(\xi)=e^{-i\frac{\xi^3}3-c\frac{\xi^2}2}.
$$
This enforces $c>0$ (if not then $\lim _{\pm \infty}\hat \psi=0$
would not hold) so that $\hat \psi$ belongs to the Schwartz space
$\mathcal S(\R)$, and so does $\psi$. The inverse Fourier
transform then yields
\begin{equation*}\label{sol-psi}
\psi(x)=\frac 1{2\pi}\int _\R \hat \psi (\xi)e^{ix\xi}\,d\xi.
\end{equation*}
But the canonical transformation yields
\begin{equation*}
  \hat \psi(\xi) = e^{-\frac{i}{3}\left(
      \xi-i\frac{c}{2}\right)^3}e^{-i\frac{c^2}{4}\xi -\frac{c^3}{24}}.
\end{equation*}
Recalling that the Airy function can be written as
\begin{equation}\label{eq:Airy}
  {\rm Ai}(x)=\frac{1}{2\pi}\int_\R e^{i\xi^3/3+ix\xi}d\xi=
  \frac{1}{2\pi}\int_{\IM \xi=\eta>0} e^{i\xi^3/3+ix\xi}d\xi ,
\end{equation}
(see e.g. \cite{Hormander}), we infer, since $\psi$ is real-valued,
\begin{align*}
  \psi(x)& = \frac 1{2\pi}\int _\R e^{-\frac{i}{3}\left(
      \xi-i\frac{c}{2}\right)^3}e^{-i\frac{c^2}{4}\xi -\frac{c^3}{24}}
  e^{ix\xi}d\xi\\
&  = \frac 1{2\pi}\int _\R e^{\frac{i}{3}\left(
      \xi+i\frac{c}{2}\right)^3}e^{i\frac{c^2}{4}\xi -\frac{c^3}{24}}
  e^{-ix\xi}d\xi\\
&= \frac 1{2\pi}\int _{\IM\zeta=c/2} e^{i\zeta^3/3}e^{i\frac{c^2}{4}(\zeta-ic/2) -\frac{c^3}{24}}
  e^{-ix(\zeta-ic/2)}d\zeta,
\end{align*}
where we have used the property $c>0$ to change the contour of
integration in the complex plane according to \eqref{eq:Airy}. Thus,
\begin{equation}\label{sol-psi3}
  \psi(x) =
e^{-cx/2+c^3/12}{\rm
Ai}\left(\frac{c^2}{4}-x\right),
\end{equation}
which is the form announced in \cite{TLK96}.

Hence, in view of \eqref{def-psi}, $\phi$ must be of the form
$\phi(x)=\psi(x-\alpha)$, with $\psi$ given by \eqref{sol-psi3}.
Conversely, if $\phi(x)=\psi(x-\alpha)$ for some $\alpha \in \R$,
then $\bar \psi =0$ enforces $ \bar \phi=\alpha$ and it is obvious
that $\phi$ solves \eqref{pb-soliton}.

\begin{theorem}[Solitary waves]\label{th:solitary waves} Let $c>0$ be given. Then there exists a
unique solitary wave $(c,\psi_c)$ solution of \eqref{pb-soliton}
and such that $\bar {\psi _c} =0$. It is given by
\eqref{sol-psi3}. Other solutions are translations of this $\psi
_c$:
$$
\phi_{\alpha,c}(x):=\psi_c(x-\alpha), \quad \alpha \in \R,
$$
so that in particular $\bar \phi_{\alpha,c}=\alpha$.

For $c\le 0$, problem \eqref{pb-soliton} has no solution.
\end{theorem}

\section{Gaussian initial data under a quadratic
  potential}
\label{sec:prop-gauss-init-quad}

Since the weight $f(x)$ in \eqref{eq-bio-poids} may be quadratic (see
e.g. \cite{LoMiPe11}), we present the explicit computations stemming
from Section~\ref{sec:algebraic-reductions} as far as the propagation
of Gaussians is concerned.
For $m\in \R$ and $a>0$, consider the Cauchy problem
\begin{equation}
  \label{eq:cauchy-quad}
  \d_t u = \d_{xx}u -\(x^2-\int_\R x^2 u(t,x)dx\)u\quad ;\quad u(0,x)
  =\sqrt{\frac{a}{2\pi}} e^{-a(x-m)^2/2}.
\end{equation}
{}From subsection \ref{subsec:gener-moment-fact},
\eqref{eq:cauchy-quad} is equivalent to
\begin{equation}
  \label{eq:cauchy-quad-red}
  \d_t v = \d_{xx}v -x^2v\quad ;\quad v(0,x) = \sqrt{\frac{a}{2\pi}} e^{-a(x-m)^2/2},
\end{equation}
through the relation
\begin{equation*}
  u(t,x) = \frac{v(t,x)}{1-\dis\int_0^t\int_\R x^2 v(s,x)dxds}.
\end{equation*}
Relation \eqref{eq:lens} shows that
\begin{equation*}
  v(t,x) = \frac{1}{\sqrt{\cosh(2t)}}
  e^{-\tanh(2t)\frac{x^2}{2}}w\(\frac{\tanh(2t)}{2},
  \frac{x}{\cosh(2t)}\),
\end{equation*}
where $w$, solution to the heat equation
\begin{equation*}
   \d_t w = \d_{xx}w\quad ;\quad w(0,x) = \sqrt{\frac{a}{2\pi}}e^{-a(x-m)^2/2},
\end{equation*}
is given by
\begin{equation*}
  w(t,x) = \sqrt{\frac{a}{2\pi(1+2at)}}e^{-\frac{a}{2(1+2at)}(x-m)^2}.
\end{equation*}
We infer
\begin{align*}
  v(t,x)& =\sqrt{ \frac{a}{2\pi(\cosh(2t) + a\sinh(2t))}}
  e^{-\tanh(2t)\frac{x^2}{2}}
  e^{-\frac{a}{2(1+a\tanh(2t))}\(\frac{x}{\cosh(2t)}-m\)^2}\\
& = \sqrt{ \frac{a}{2\pi(\cosh(2t) + a\sinh(2t))}}
  e^{-\frac{am^2 \sinh(2t)}{2(a\cosh(2t)+\sinh(2t))}}
  e^{-\frac{a(t)}{2}
      \(x-m(t)\)^2},
\end{align*}
where
\begin{equation*}
  a(t) = \frac{a\cosh(2t)+\sinh(2t)}{\cosh(2t)+a\sinh(2t)},\quad m(t)
  = \frac{am}{a\cosh(2t)+\sinh(2t)},
\end{equation*}
hence
\begin{multline*}
  \int_\R x^2v(t,x)dx = \frac{\sqrt a e^{-\frac{am^2
        \sinh(2t)}{2(a\cosh(2t)+\sinh(2t))}}}{\(a\cosh(2t) +
      \sinh(2t)\)^{5/2} }\\
     \times  \((\cosh(2t)+a\sinh(2t))(a\cosh(2t)+\sinh(2t))+a^2m^2\).
\end{multline*}
The integral in time of this quantity involves elliptic integrals in
general, so we consider special values of the parameters.
In the particular case of an initial Gaussian centered at the origin, $m=0$,
with $a=1$, the above formula becomes much simpler,
\begin{equation*}
   \int_\R x^2v(t,x)dx =  e^{-t},
\end{equation*}
and we check that
\begin{equation}\label{eq:ground}
  u(t,x) = \frac{1}{\sqrt{2\pi}}e^{-x^2/2}= u_0(x),\quad \forall t\in
  \R.
\end{equation}
In other cases, the initial Gaussian propagates as a Gaussian in a
non-trivial way, for which in general explicit computations seem
rather intricate. The fact that the solution in \eqref{eq:ground} does
not depend on time can be understood as follows: the Gaussian
$e^{-x^2/2}$ is the ground state associated to the harmonic
oscillator, that is the eigenfunction associated to the lowest
eigenvalue of the harmonic oscillator (see e.g. \cite{LandauQ}):
\begin{equation*}
  \(-\d_{xx}+x^2\) e^{-x^2/2} = e^{-x^2/2} ,
\end{equation*}
so the solution to \eqref{eq:cauchy-quad-red} is simply
\begin{equation*}
  v(t,x)=\frac{1}{\sqrt{2\pi}}e^{-t}e^{-x^2/2},
\end{equation*}
hence
\begin{equation*}
    \int_\R x^2v(t,x)dx =  e^{-t},
\end{equation*}
and $u(t,x)=u_0(x)$ from \eqref{eq:chgt1}. Note that this specific
case (stationary solution)
does not extend to other biologically relevant cases: the
eigenfunctions associated to the harmonic oscillator are Hermite
functions,
\begin{equation*}
  \psi_n(x) = e^{x^2/2}\frac{d^n}{dx^n}e^{-x^2},
\end{equation*}
which are associated with the eigenvalue $\lambda_n=1+2n$, but except
in the case $n=0$, they change signs.
\bigskip

\noindent {\bf Acknowledgements.} The authors are grateful to
J\'er\^ome Coville for pointing out the models studied in the
present paper. They also would like to thank Sylvain Gandon and
Ga\"el Raoul for valuable discussions, and for bringing reference
\cite{B14} to their attention. M. A. is supported by the French
{\it Agence Nationale de la Recherche} within the project IDEE
(ANR-2010-0112-01).

\bibliographystyle{siam}
\bibliography{biblio}

\end{document}